\documentclass[12pt]{article}
\usepackage{bbm}
\usepackage{amsfonts}
\usepackage{pifont}
\usepackage{hyperref}
\usepackage{mathrsfs}
\usepackage{indentfirst}
\usepackage{amsmath}
\usepackage{amssymb}
\usepackage[amsmath, thmmarks]{ntheorem}
\usepackage{cite}
\usepackage{graphicx}
\usepackage{tikz}
\newtheorem{definition}{\bf Definition}[section]
\newtheorem{lemma}{\bf Lemma}[section]
\newtheorem{theorem}{\bf Theorem}[section]
\newtheorem{remark}{\bf Remark}[section]
\newtheorem{corollary}{\bf Corollary}[section]
\newtheorem{proposition}{\bf Proposition}[section]
\newtheorem{example}{\bf Example}[section]

\def\QEDopen{{\setlength{\fboxsep}{0pt}\setlength{\fboxrule}{0.2pt}\fbox{\rule[0pt]{0pt}{1.3ex}\rule[0pt]{1.3ex}{0pt}}}} %¶¨Òå¿ÕÐÄ·û
\def\QED{\QEDopen}
\def\proof{{\bf Proof.} }
\def\endproof{\hspace*{\fill}~\QED\par\endtrivlist\unskip}

\begin{document}
\setcounter{page}{1}

\title{{\textbf{Structures of lattices which can be represented as the collection of all up-sets}}\thanks {Supported by National Natural Science
Foundation of China (No.61573240)}}
\author{Peng He\footnote{\emph{E-mail address}: 443966297@qq.com}, Xue-ping Wang\footnote{Corresponding author. xpwang1@hotmail.com; fax: +86-28-84761393}\\
\emph{College of Mathematics and Software
Science, Sichuan Normal University,}\\
\emph{Chengdu, Sichuan 610066, People's Republic of China}}

\newcommand{\pp}[2]{\frac{\partial #1}{\partial #2}}
\date{}
\maketitle
\begin{quote}
{\bf Abstract}
This paper first gives a necessary and sufficient condition that a lattice $L$ can be represented as the collection of all up-sets of a poset. Applying the condition, it obtains a necessary and sufficient condition that a lattice can be embedded into the lattice $L$ such that all infima, suprema, the top and bottom elements are preserved under the embedding by defining a monotonic operator on a poset. This paper finally shows that the quotient of the set of the monotonic operators under an equivalence relation can be naturally ordered and it is a lattice if $L$ is a finite distributive lattice.


\end{quote}
\section{Introduction}\label{intro}

$L$-fuzzy sets and structures have been widely studied from Goguen's first paper \cite{Goguen}. These structures appear
when the membership grades can be represented by elements of an ordered set, instead of just by numbers in the unit $[0,1]$.
$L$-fuzzy sets play an important role in many areas of research such as algebraic theories including order-theoretic
structures (see e.g., \cite{Seselja5}), automata and tree series (see e.g., \cite{Borchar}) and theoretical computer
science (see e.g., \cite{Chechik}). It is well known that every $L$-fuzzy structure is uniquely determined by the collection
of cut sets. The cut sets can be ordered naturally by reverse inclusion. Thus cut sets are one of the most important
links between $L$-fuzzy mathematics and classical theory of ordered structures. That is to say,
$L$-fuzzy mathematics provides many useful techniques and methods to investigate classical theory of ordered structures (see e.g., \cite{Seselja2, Seselja98}).

In 2003, \v{S}e\v{s}elja and Tepav\v{c}evi\'{c} \cite{Seselja3} introduced a particular completion which is equivalent with the famous
Dedekind-MacNeille completion by fuzzy sets. Moreover, they presented a survey on representations of ordered
structures by fuzzy sets, and proved that the structure itself is uniquely represented by the collection of cut sets
ordered dually to inclusion (see \cite{Seselja4}).

In 2010, Jim\'{e}nez, Montes, \v{S}e\v{s}elja and Tepav\v{c}evi\'{c} \cite{Jorge} first introduced a $L$-fuzzy up-set and then they
showed a necessary and sufficient condition under which a collection of crisp up-sets of a poset $X$ consists of cut sets of the $L$-fuzzy up-set. In particular, they obtained a representation of any finite distributive lattice as a family of cut sets of an $L$-fuzzy up-set, i.e., they gave a version
of the famous Bikhoff Representation Theorem by using $L$-fuzzy up-sets. The famous Bikhoff Representation Theorem says that any finite distributive lattice $L$ can be isomorphically represented by the collection of all up-sets on the set of all meet-irreducible elements of $L$. Therefore, an interesting problem is: What is the structure of a lattice which can be represented as the collection of all up-sets on the set of all meet-irreducible elements of $L$? This paper will focus on the problem by applying $L$-fuzzy up-sets.

The paper is organized as follows. For the sake of convenience, it gives some notions and previous results in Section 2.
In Section 3, it obtains a necessary and sufficient condition under which a lattice can be represented as the collection of all up-sets on the set of all completely meet-irreducible elements of $L$ by using $L$-fuzzy up-sets.
In Section 4, it first introduces a concept of a monotonic operator on a poset,
and then presents a necessary and sufficient condition that a lattice can be embedded into a given lattice which can be represented as the collection of all up-sets. This paper finally shows that the quotient of the monotonic operators on a finite distributive lattice forms a lattice in Section 5. Conclusions are drawn in Section 6.

\section{Preliminaries}
We first list some necessary notions and relevant properties from the classical order theory in the sequel.
For more comprehensive presentation, see e.g., books \cite{Crawley, Gratze}.

A poset is a structure $(P, \leq)$, or $P$ for short, where $P$ is a nonempty set and $\leq$ an ordering (reflexive, antisymmetric and transitive)
relation on $P$. A complete lattice is a poset $(L, \leq)$ in which every subset $M$ has the greatest lower bound, infimum,
meet, denoted by $\bigwedge M$, and the least upper bound, supremum, join, denoted by $\bigvee M$. A complete lattice $L$ possesses the top element
$1_L$ and the bottom element $0_L$. We say that an element $a$ in a lattice $L$ is completely meet-irreducible if $a\neq 1_L$
and from every family $\{x_i\mid i\in I\}$ of elements in $L$, from $a=\bigwedge_{i\in I}x_i$ it follows that $a=x_i$ for some $i\in I$.
Furthermore, we denote by $M(L)$ the set of all completely meet-irreducible elements of $L$.

A complete lattice $L$ is called infinitely distributive if, for all $x\in L$ and $S\subseteq L$, $x\wedge \bigvee S=\bigvee_{s\in S}x\wedge s$.
A complete lattice $L$ is called dual infinitely distributive if its dual is infinitely distributive.

An up-set (a semi-filter) on a poset $P$
is any sub-poset $U$, satisfying: for $x\in U, y\in P, x\leq y$ implies $y\in U$.

Next, we recall a natural way that distributive lattices appear
among ordered structures.

\begin{lemma}\label{lemm 0}
The collection of all up-sets of a poset $X$ is a complete distributive lattice under inclusion. Further, it is infinitely and dual infinitely distributive.
\end{lemma}

If $a$ is an element of a complete lattice $L$, then a representation $a=\bigwedge Q$ with $Q\subseteq M(L)$ is called
a decomposition of $a$. Further, we say a complete lattice $L$ has a \emph{Decomposition Property} (often abbreviated as $DP$)
if every element of $L$ has a decomposition.

In the following, we present some notations from theory of $L$-fuzzy up-sets. More details about the relevant properties can be found e.g.,
in \cite{Jorge, Seselja3, Seselja4}.

An $L$-fuzzy up-set is a mapping $\mu : X\rightarrow L$ from a poset $(X, \leq)$ (domain) into a complete
lattice $(L, \leq)$ (co-domain) with the condition that for all $x, y\in X$
$$x\leq y \mbox{ implies } \mu(x)\leq \mu(y).$$

If $\mu : X\rightarrow L$ is an $L$-fuzzy up-set on poset $X$ then, for $p\in L$, the set \begin{equation*}\label{equ 00}\mu_p =\{x\in X\mid \mu(x)\geq p\}\end{equation*}
is called the \emph{$p$-cut}, a \emph{cut set} or simply a \emph{cut} of $\mu$. Let $\mu_L=\{\mu_p\mid p\in L\}$. Moreover, let $L^{\mu}$ be
defined by \begin{equation}\label{equ 0}L^{\mu}=\{p\in L\mid p=\bigwedge B \mbox{ with }
B\subseteq \mu(X)\},\end{equation} where $\mu(X)=\{\mu(x)\mid x\in X\}$. Then $L^{\mu}$ is a complete lattice (see \cite{Jorge}).

The following four statements present some characterizations of the collection of cuts of an $L$-fuzzy up-set.
\begin{proposition}[\cite{Jorge, Seselja3}]\label{pro 1}
Let $\mathcal{F}$ be a family of some up-sets of a poset $X$ which is closed under intersections and contains $X$ and $F=(\mathcal{F}, \supseteq)$. Let $\mu: X\rightarrow F$ be defined by
$$\mu(x)=\bigcap \{p\in\mathcal{F}\mid x\in p\}.$$ Then, $\mu$ is an $F$-fuzzy up-set on $X$ with the co-domain lattice $F$ such that
$\mu_F=\mathcal{F}$ and for every $p\in \mathcal{F}$ it holds that $p=\mu_p$.
\end{proposition}
\begin{lemma}[\cite{He}]\label{lemm 1}
Let $\mathcal{F}$ be a family of some up-sets of a poset $X$ which is closed under intersections and contains $X$ and let $F=(\mathcal{F}, \supseteq)$.
Then there exists an $F$-fuzzy up-set $\mu: X\rightarrow F$ such that $\mu_F= \mathcal{F}$ and $F^{\mu}=F$.
\end{lemma}
\begin{proposition}[\cite{Jorge}]\label{pro 2}
Let $X$ be a poset, let $L$ be a complete lattice and let $\mu: X\rightarrow L$ be an $L$-fuzzy up-set. Then the following two statements are equivalent: \\
\emph{(a1)} $\mu_L$ is formed by all the up-sets of $X$.\\
\emph{(a2)} For every family $\{x_i\mid i\in I\}$ of elements from $X$ and every $x\in X$, it holds that
$$\mu(x)\geq \bigwedge_{i\in I} \mu(x_i) \mbox{ implies }x\geq x_i \mbox{ for some } i\in I.$$
\end{proposition}
\begin{proposition}[\cite{Jorge}]\label{pro 3}
Let $X$ be a poset, let $L$ be a complete lattice and let $\mu$: $X\rightarrow L$ be an $L$-fuzzy up-set. If $\mu_L$ is
formed by all the up-sets of $X$, then $\mu(x)\in M(L^{\mu})$ for all $x\in X$.
\end{proposition}

Given any set $A$, we denote its cardinality by $|A|$. Then we say a lattice $L$ is trivial if $|L|=1$.
In what follows, we assume that all lattices are non-trivial and denote by $\mathcal{L}_{CD}$ the set of all complete distributive lattices.

\section{Representation of a complete distributive lattice by $L$-fuzzy up-sets}
In this section, we shall show a necessary and sufficient condition under which a lattice can be represented as the collection of all up-sets on the set of all completely meet-irreducible elements of $L$ by using $L$-fuzzy up-sets.

For convenience, let $\mathcal{F}_X$ be the family of all up-sets of a poset $X$ and denote $(\mathcal{F}_X, \supseteq)$ by $F_X$. We first give the following two lemmas.

\begin{lemma}\label{lemm 2}
Let $X$ be a poset. Then $F_X\in \mathcal{L}_{CD}$ and $F_X$ has $DP$ and satisfies the following condition:

\emph{($\mathfrak{M}$)} For each $q\in M(F_X)$, if $q\geq \bigwedge_{C\subseteq M(F_X)} C$ then $q\geq p$ for some $p\in C$.
\end{lemma}
\proof
By Lemma \ref{lemm 0}, $F_X\in \mathcal{L}_{CD}$. By Lemma \ref{lemm 1},
we have an $F_X$-fuzzy up-set $$\mu: X\rightarrow F_X$$ such that $F_X^{\mu}=F_X$ and
$\mu_{F_X}=\mathcal{F}_X$. Further, by Proposition \ref{pro 3}, $\mu(x)\in M(F_X^{\mu})$ for all $x\in X$.
Then it follows from formula (\ref{equ 0}) that $F_X^{\mu}$ has $DP$, which means that $F_X$ has $DP$.

Next, we shall prove that the condition ($\mathfrak{M}$) holds.
Let $q\in M(F_X)$ and $q\geq \bigwedge_{C\subseteq M(F_X)} C$. Then $$q=q\vee\bigwedge_{C\subseteq M(F_X)} C=\bigwedge_{p\in C} (q\vee p)$$
by Lemma \ref{lemm 0}. Thus, $q=q\vee p$ for some $p\in C$ since $q\in M(F_X)$. Therefore, $q\geq p$, i.e., ($\mathfrak{M}$) holds.
\endproof

\begin{lemma}\label{lemm 3}
Let $L\in \mathcal{L}_{CD}$ and $L$ have $DP$ and satisfy \emph{($\mathfrak{M}$)}. Then $L\cong F_{M(L)}$.
\end{lemma}
\proof
Define $\mu: M(L)\rightarrow L$ by $\mu(x)=x$ for all $x\in M(L)$. Clearly, $\mu$ is an $L$-fuzzy up-set. Suppose that
$\mu(x)\geq \bigwedge_{i\in I}\mu(x_i)$. Then $x\geq \bigwedge_{i\in I}x_i$.
Thus the condition ($\mathfrak{M}$) yields that there exists $i\in I$ such that $x\geq x_i$.
Further, from Proposition \ref{pro 2}, it follows that $\mathcal{F}_{M(L)}=\mu_L$. Therefore,
it suffices to prove $L\cong(\mu_L, \supseteq)$. Let $f: L\rightarrow (\mu_L, \supseteq)$ be defined by $f(p)=\mu_p$
for all $p\in L$. We just need to prove that $f$ is isomorphic.

First, $f$ is a map from $L$ to $(\mu_L, \supseteq)$ obviously. Furthermore, if $T\in \mu_L$ then $T=\mu_q$ for some $q\in L$.
Thus $T =\mu_q =f(q)$, which means that $f$ is surjective.

Secondly, assume that $p\neq q$ but $f(p)=f(q)$.
Let $r\in L$. We have that \begin{equation}\label{equ 4} f(r)=\mu_r=\{x\in M(L)\mid x\geq r\}\end{equation} since $\mu(x)=x$ for all $x\in M(L)$.
Furthermore, as $L$ has $DP$, \begin{equation}\label{equ hhh4} r=\bigwedge f(r).\end{equation}
Thus $f(p)=f(q)$ implies that $p=q$, a contradiction. Therefore, $f$ is injective.

In view of the proof in the preceding two paragraphs, the mapping $f$ is a bijection.

Finally, we prove that, $f$ and $f^{-1}$, preserve the orders $\leq$ and $\supseteq$.

If $p\leq q$ then obviously $\mu_p\supseteq \mu_q$, i.e.,
$f(p)\supseteq f(q)$, which means that $f$ preserves the order $\leq$.

Let $S\in \mu_L$. Then $S=\mu_s$ for some $s\in L$. Thus $S=f(s)$  and $s=\bigwedge S$ by formulas (\ref{equ 4}) and (\ref{equ hhh4}).
So that $f^{-1}(S)=\bigwedge S$ for each $S\in \mu_L$. If $T_1\subseteq T_2$ and $T_1, T_2 \in \mu_L$ then $f^{-1}(T_1)=\bigwedge T_1 \geq \bigwedge T_2 = f^{-1}(T_2)$.
Therefore, $f^{-1}$ also preserves the order $\supseteq$.

Thus $f$ is an isomorphism from $L$ to $(\mu_L, \supseteq)$. Consequently, $L\cong F_{M(L)}$.
\endproof

The following corollary is directly by Lemma \ref{lemm 3}.

\begin{corollary}\label{cor3}
Let $L_1, L_2\in \mathcal{L}_{CD}$ and both of them have $DP$ and satisfy \emph{($\mathfrak{M}$)}. Then $L_1 \cong L_2$ if and only if $(M(L_1), \leq)\cong (M(L_2), \leq)$.
\end{corollary}

\begin{theorem}\label{theo 5}
The following conditions are equivalent:\\
\emph{(b1)} There exists a post $X$ such that $L\cong F_X$.\\
\emph{(b2)} $L\in \mathcal{L}_{CD}$ and $L$ has $DP$ and satisfies \emph{($\mathfrak{M}$)}.\\
\emph{(b3)} $L$ is dual infinitely distributive and it has $DP$.
\end{theorem}
\proof
First, it is clear that (b1) and (b2) are equivalent by Lemmas \ref{lemm 2} and \ref{lemm 3}.
Secondly, it follows from Lemmas \ref{lemm 0} and \ref{lemm 2} that (b1) implies (b3).
Finally, we shall prove that (b3) implies (b2). Obviously, the condition that $L$ is dual infinitely distributive yields that
$L\in \mathcal{L}_{CD}$. It remains, therefore, to show that (b3) implies ($\mathfrak{M}$).
Let $q\in M(L)$ and $\{q_i\mid i\in I\}\subseteq M(L)$. Suppose that $q\geq \bigwedge_{i\in I}q_{i}$.
Then $$q=q\vee \bigwedge_{i\in I}q_{i}=\bigwedge_{i\in I} q\vee q_i$$ since $L$ is dual infinitely distributive.
Further, as $q\in M(L)$, $q=q\vee q_j$ for some $j\in I$. Then $q\geq q_j$. Thus, the condition ($\mathfrak{M}$) holds.
Consequently, (b3) implies (b2).
\endproof

By Lemma \ref{lemm 0} and Theorem \ref{theo 5}, we have
\begin{corollary}\label{corx}
If a lattice $L$ is dual infinitely distributive and it has $DP$, then it is also infinitely distributive.
\end{corollary}

Note that not every dual infinitely distributive lattice is infinitely distributive (see p.118 in \cite{Birkhoff}). Also, not every infinitely distributive and dual infinitely distributive lattice has $DP$.
For example, the complete chain $[0, 1]$ does not have $DP$.

Applying Theorem \ref{theo 5} again, we have
\begin{corollary}\label{cor4}
Let $L\in \mathcal{L}_{CD}$. If $L$ is finite then $L$ has $DP$ and satisfy \emph{($\mathfrak{M}$)}.
\end{corollary}

\begin{proposition}\label{propo1}
Let $L\in \mathcal{L}_{CD}$ and $L$ have $DP$ and satisfy \emph{($\mathfrak{M}$)}. If $|M(L)|=n$ then
each maximal chain of $L$ has $n+1$ elements.
\end{proposition}
\proof
Let $M(L)=\{a_1, a_2, \cdots, a_n\}$. Without loss of generality, suppose that \begin{equation}\label{he4441}a_i\ngtr a_j\end{equation} for
any $1\leq i< j\leq n$. Let $C_k=\{a_k, \cdots, a_n\}$ if $1\leq k\leq n$, otherwise, let $C_{n+1}=\emptyset$.

We claim that $C_k$ is an up-set on $(M(L), \leq)$ for any $1\leq k\leq n+1$.
First, $C_{n+1}$ is an up-set. Now, we prove that $C_k$ is an up-set for any $1\leq k\leq n$.
Let $a_u\in M(L)$, $a_v\in C_k$ and $a_u\geq a_v$. By formula (\ref{he4441}), $u\geq v$.
On the other hand, $a_v\in C_k$ deduces that $v\geq k$. Thus $u\geq k$, and then $a_u\in C_k$.
So, $C_k$ is an up-set. Therefore, $C_k$ is an up-set on $(M(L), \leq)$ for any $1\leq k\leq n+1$.

By Lemma \ref{lemm 3}, $L\cong F_{M(L)}$. Obviously, $M(L)=C_1\prec C_2\prec  \cdots\prec  C_{n+1}=\emptyset$ is a
maximal chain of $F_{M(L)}$. So that there exists a maximal chain of $L$ which has $n+1$ elements.
As $L\in \mathcal{L}_{CD}$, $|l|=n+1$ for each maximal chain $l$ of $L$.
\endproof

From Proposition \ref{propo1}, we have

\begin{remark}\label{remarkk1}
\emph{If $L$ is an infinite complete distributive lattice and
it has $DP$ and satisfies ($\mathfrak{M}$), then there exists an infinite chain in $L$.}
\end{remark}

Let us conclude this section with a version of the famous Birkhoff Representation Theorem.
From Theorem \ref{theo 5}, we know that any complete distributive lattice $L$ can be isomorphically represented by the collection
of all up-sets on $(M(L), \leq)$ if and only if $L$ has $DP$ and satisfies ($\mathfrak{M}$).
Moreover, by the proof of Lemma \ref{lemm 3}, we can define an $L$-fuzzy up-set whose cut sets are order isomorphic to $L$.
Therefore, we obtain a representation of the lattice $L$ as the family of all cut sets of an $L$-fuzzy up-set. This $L$-fuzzy up-set is
the embedding of $M(L)$ in $L$, that is, $\mu: M(L)\rightarrow L$ with $\mu(x)=x$ for every $x\in M(L)$.

\section{Embeddedness of a class of distributive lattices}
In what follows, let $E(L)$ be the set of all lattices which can be embedded into the lattice $L$, such that
all infima, suprema, the top and bottom elements are preserved under the embedding.

We know that every complete sublattice of a complete distributive lattice $L$ can be embedded into $L$ such
that all infima and superma are preserved under the embedding. However, the sublattice may not be in $E(L)$ as illustrated by the following example.
\begin{example}\label{exa 1}
\emph{ Let us consider the complete distributive lattices $L_0$ and $L$ represented in Fig. 1.}
\par\noindent\vskip50pt
 \begin{minipage}{11pc}
\setlength{\unitlength}{0.75pt}\begin{picture}(600,100)
\put(330,40){\circle{4}}\put(335,40){\makebox(0,0)[l]{\footnotesize $0_L$}}
\put(330,75){\circle{4}}\put(335,65){\makebox(0,0)[l]{\footnotesize $d$}}
\put(295,110){\circle{4}}\put(290,120){\makebox(0,0)[l]{\footnotesize $a$}}
\put(365,110){\circle{4}}\put(370,120){\makebox(0,0)[l]{\footnotesize $c$}}
\put(330,145){\circle{4}}\put(335,150){\makebox(0,0)[l]{\footnotesize $1_L$}}
\put(330,42){\line(0,1){31}}
\put(331.5,76.5){\line(1,1){32}}
\put(328.5,76.5){\line(-1,1){32}}
\put(296.5,111.5){\line(1,1){32}}
\put(363.5,111.5){\line(-1,1){32}}
\put(325,20){$L$}
\put(165,75){\circle{4}}\put(160,65){\makebox(0,0)[l]{\footnotesize $a$}}
\put(235,75){\circle{4}}\put(240,65){\makebox(0,0)[l]{\footnotesize $c$}}
\put(200,110){\circle{4}}\put(205,115){\makebox(0,0)[l]{\footnotesize $1_L$}}
\put(200,40){\circle{4}}\put(205,38){\makebox(0,0)[l]{\footnotesize $d$}}
\put(201.5,41.5){\line(1,1){32}}
\put(198.5,41.5){\line(-1,1){32}}
\put(166.5,76.5){\line(1,1){32}}
\put(233.5,76.5){\line(-1,1){32}}
\put(190,20){$L_0$}
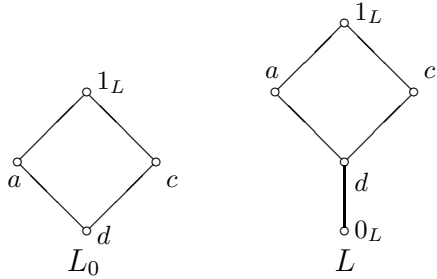
\put(155,0){\emph{Fig.1 Hasse diagrams of $L_0$ and $L$}}
  \end{picture}
  \end{minipage}\\
\emph{Clearly, $L_0$ is a sublattice of $L$ but $L_0\notin E(L)$.}
\end{example}

Let $L_0, L\in \mathcal{L}_{CD}$ and they have $DP$ and satisfy ($\mathfrak{M}$). Then from Lemma \ref{lemm 3}, $L\cong F_{M(L)}$.
In \cite{He1}, He and Wang proved that for all closure operators $C$ on $M(L)$, if $L_0\cong F_{M(L)/C}$ then $L_0 \in E(L)$.
Also, they gave an example to show that there exists $L_1 \in E(L)$ such that $L_1\ncong F_{M(L)/C_1}$ for any
closure operators $C_1$ on $M(L)$, in which $L_1\in \mathcal{L}_{CD}$ and $L_1$ has $DP$ and satisfies ($\mathfrak{M}$). Therefore, an interesting problem is: What is the condition that $L_1\cong F_{M(L)/C^*}$ for a certain operator $C^*$ on $M(L)$ for any $L_1 \in E(L)$?

In this section, we shall define a monotonic operator $G$ on $M(L)$, and then prove that $L_0\in E(L)$
if and only if there exists a monotonic operator $G$ on $M(L)$ such that $L_0\cong F_{M(L)/G}$.

Let $X$ be a nonempty set and $\mathcal{P}(X)=\{S\mid S\subseteq X\}$. Furthermore, we denote $(\mathcal{P}(X), \supseteq)$ by $P_X$.
\begin{definition}\label{def 1}
\emph{A monotonic operator $G$ on poset $X$ is a function $G: X\rightarrow P_X$ such that, for all $x, y \in X$, $x\leq y$ implies $G(x)\supseteq G(y)$.}
\end{definition}

Clearly, a monotonic operator $G$ is also a $P_X$-fuzzy up-set on $X$ and we have the following lemma.
\begin{lemma}\label{lemma 1}
Let $G$ be a monotonic operator on poset $X$. We consider a relation $\sim$ on $X$ defined by $x\sim y$ if and only if $G(x)= G(y)$. Then:\\
\emph{(a)} $\sim$ is an equivalence relation on $X$.\\
\emph{(b)} The set $X/\sim$ can be ordered: $[x]\leq [y]$ if and only if $G(x)\supseteq G(y)$ where $[x]=\{z\in X\mid x\sim z\}$ for all $x\in X$.\\
\emph{(c)} $(X/\sim, \leq)\cong(G(X), \supseteq)$ where $G(X)=\{G(x)\mid x\in X\}$.
\end{lemma}
\proof
Obviously, (a) and (b) hold. Now, we shall show (c).

Let $f: (X/\sim, \leq)\rightarrow (G(X), \supseteq)$ be defined by
$$f([x])=G(x)$$ for all $[x]\in X/\sim$. Thus, we only need to prove $f$ is an isomorphism.
From (b), $[x]=[y]$ if and only if $G(x)=G(y)$, then $f$ is a map from $(X/\sim, \leq)$ to $(G(X), \supseteq)$.
Further, $f$ is bijective obviously. Again, from (b), we know that both $f$ and $f^{-1}$ preserve the orders $\leq$ and $\supseteq$, respectively.
Therefore, $f$ is isomorphic.
\endproof

In what follows, we denote $X/\sim$ by $X/G$ if $G$ is a monotonic operator on $X$. Then we have the theorem as below.

\begin{theorem}\label{the0 6}
Let $X$ be a poset and $G$ be a monotonic operator on $X$. Then $F_{X/G} \in E(F_X)$.
\end{theorem}
\proof
Let \begin{equation}\label{equa 4}S_T=\bigcup_{[x]\in T}\{y\in X\mid y\in [x]\}\end{equation} for all $T\in \mathcal{F}_{X/G}$ and
\begin{equation}\label{equa 5}\mathcal{S}_G=\{S_T\mid T\in \mathcal{F}_{X/G}\}.\end{equation}
Now, we prove $F_{X/G} \cong (\mathcal{S}_G, \supseteq)$. Let $f: F_{X/G}\rightarrow (\mathcal{S}_G, \supseteq)$ be defined
by $f(T)=S_T$ for all $T\in \mathcal{F}_{X/G}$.
Clearly, $f$ is a surjective map by (\ref{equa 4}) and (\ref{equa 5}). Suppose that $T_1 \neq T_2$ in $\mathcal{S}_G$. We claim that $f(T_1) \neq f(T_2)$.
Otherwise, $f(T_1)=f(T_2)$. Then $S_{T_1}= S_{T_2}$. Without loss of generality, we suppose that $[x]\in T_1$ but
$[x]\notin T_2$ since $T_1 \neq T_2$. Thus $x\in [x] \subseteq S_{T_1}$ by (\ref{equa 4}), which means that $x\in S_{T_2}$. Hence, there exists $[y]\in T_2$
such that $x\in [y]$, and this means that $[x]=[y]\in T_2$, a contradiction. So that $f$ is injective.

Therefore, $f$ is bijective.

We prove that both $f$ and $f^{-1}$ preserve the order $\supseteq$.

If $T_1 \supseteq T_2$ then obviously $f(T_1)=S_{T_1}\supseteq S_{T_2}=f(T_2)$, i.e., $f$ preserves the order $\supseteq$.

Furthermore, we claim that $T_1 \supseteq T_2$ when $f(T_1)\supseteq f(T_2)$. Suppose $T_1 \nsupseteq T_2$ but $f(T_1)\supseteq f(T_2)$.
Then there exists an element $[z]\in T_2$ but $[z]\notin T_1$. Note that $z\in S_{T_2}=f(T_2)$ by formula (\ref{equa 4}). This implies that
$z\in f(T_1)$ since $f(T_1)\supseteq f(T_2)$. Thus, from (\ref{equa 4}), there exists $[u]\in T_1$ such
that $z\in [u]$, and so that $[z]=[u]\in T_1$, a contradiction. Consequently, $f^{-1}$ preserves the order $\supseteq$.

Thus $f$ is an isomorphism from $F_{X/G}$ to $(\mathcal{S}_G, \supseteq)$, i.e., \begin{equation}\label{equa 7}F_{X/G} \cong (\mathcal{S}_G, \supseteq).\end{equation}

In what follows, we shall prove that $F_{X/G}\in E(F_X)$.

First, by (\ref{equa 4}), we observe that \begin{equation}\label{eq112}\emptyset=S_{\emptyset}\in \mathcal{S}_G\mbox{ and }X=S_{X/G}\in \mathcal{S}_G.\end{equation} Now,
let $K\in \mathcal{S}_G$. Then, applying (\ref{equa 5}), we know that there is a $T \in \mathcal{F}_{X/G}$ such that $K=S_T$. By Lemma \ref{lemma 1} and Definition \ref{def 1}, we
know that for all $x, y\in X$, \begin{equation}\label{2222}y\geq x \mbox{ in }X\mbox{ implies } [y]\geq [x] \mbox{ in } X/G.\end{equation}
Suppose that $x\in K$, $y\in X$ and $y\geq x$. Then from formula (\ref{equa 4}), we have $[x]\in T$ since $K=S_T$.
Thus, by formula (\ref{2222}), the condition $y\geq x$ implies that $[y]\in T$ since $T\in \mathcal{F}_{X/G}$. Hence, $y\in S_T=K$ which means that $K\in \mathcal{F}_X$.
Therefore, by the arbitrariness of $K$, we
have \begin{equation}\label{equa 9}\mathcal{S}_G\subseteq \mathcal{F}_X.\end{equation}

Secondly, let $\{S_i\mid i\in I\}\subseteq\mathcal{S}_G$. From (\ref{equa 5}), it follows that there
exists $T_i \in \mathcal{F}_{X/G}$ such that $S_i=S_{T_{i}}$ for any $i\in I$. Thus $\bigcap_{i\in I}S_i=\bigcap_{i\in I}S_{T_{i}}$
and $\bigcup_{i\in I}S_i=\bigcup_{i\in I}S_{T_{i}}$. Let $W=\bigcap_{i\in I}T_i$ and $R=\bigcup_{i\in I}T_i$. Observe that $W, R\in \mathcal{F}_{X/G}$
since $F_{X/G}$ is a complete lattice. Thus, by (\ref{equa 5}), $S_W, S_R \in \mathcal{S}_G$.
Note that $[x]\cap [y]=\emptyset$ for any $[x]\neq [y]$. Thus, by (\ref{equa 4}), we further know that $S_W=\bigcap_{i\in I}S_{T_i}$ and $S_R=\bigcup_{i\in I}S_{T_i}$.
Therefore, \begin{equation}\label{equa 8}\bigcap_{i\in I}S_i \in \mathcal{S}_G \mbox{ and } \bigcup_{i\in I}S_i \in \mathcal{S}_G.\end{equation}

Finally, from formulas (\ref{equa 7}), (\ref{eq112}), (\ref{equa 9}) and (\ref{equa 8}), we know that $F_{X/G}\in E(F_X)$.
\endproof

The following example will illustrate Theorem \ref{the0 6}.

\begin{example}\label{exam 2}
\emph{Let us consider the posts $X$ and $X/G$ represented in Fig. 2, where $G$ is a is a monotonic operator on $X$ with
$G(x)=\left\{
\begin{array}{rcl}
a_1     &      & {x=a,}\\
b_1     &      & {x=b,}\\
c_1     &      & {x=c.}
\end{array}\right.$ }
\par\noindent\vskip50pt
 \begin{minipage}{11pc}
\setlength{\unitlength}{0.75pt}\begin{picture}(600,100)
\put(280,40){\circle{4}}\put(285,38){\makebox(0,0)[l]{\footnotesize $a_2$}}
\put(245,75){\circle{4}}\put(240,65){\makebox(0,0)[l]{\footnotesize $a_1$}}
\put(280,75){\circle{4}}\put(285,73){\makebox(0,0)[l]{\footnotesize $c_2$}}
\put(315,75){\circle{4}}
\put(280,110){\circle{4}}
\put(245,110){\circle{4}}\put(235,120){\makebox(0,0)[l]{\footnotesize $c_1$}}
\put(315,110){\circle{4}}\put(320,120){\makebox(0,0)[l]{\footnotesize $b_2$}}
\put(280,145){\circle{4}}\put(283,150){\makebox(0,0)[l]{\footnotesize $b_1$}}
\put(245,77){\line(0,1){31}}
\put(280,112){\line(0,1){31}}
\put(278.5,41.5){\line(-1,1){32}}
\put(281.5,41.5){\line(1,1){32}}
\put(280,42){\line(0,1){31}}
\put(281.5,76.5){\line(1,1){32}}
\put(278.5,76.5){\line(-1,1){32}}
\put(315,77){\line(0,1){31}}
\put(313.5,76.5){\line(-1,1){32}}
\put(246.5,76.5){\line(1,1){32}}
\put(246.5,111.5){\line(1,1){32}}
\put(313.5,111.5){\line(-1,1){32}}
\put(270,20){$P_X$}
\put(115,40){\circle{4}}\put(110,30){\makebox(0,0)[l]{\footnotesize $a$}}
\put(185,40){\circle{4}}\put(190,30){\makebox(0,0)[l]{\footnotesize $c$}}
\put(150,75){\circle{4}}\put(155,80){\makebox(0,0)[l]{\footnotesize $b$}}
\put(116.5,41.5){\line(1,1){32}}
\put(183.5,41.5){\line(-1,1){32}}
\put(144,20){$X$}
\put(370,110){\circle{4}}\put(373,115){\makebox(0,0)[l]{\footnotesize $[b]$}}
\put(370,75){\circle{4}}\put(352,80){\makebox(0,0)[l]{\footnotesize $[c]$}}
\put(370,40){\circle{4}}\put(373,45){\makebox(0,0)[l]{\footnotesize $[a]$}}
\put(370,42){\line(0,1){31}}
\put(370,77){\line(0,1){31}}
\put(356,20){$X/G$}
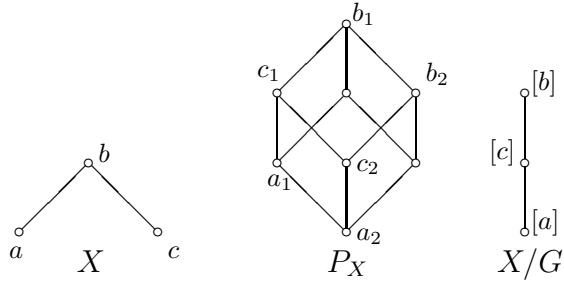
\put(111,0){\emph{Fig.2 Hasse diagrams of $X$, $P_X$ and $X/G$}}
  \end{picture}
  \end{minipage}\\
\emph{It is easy to check that $$\mathcal{F}_X=\{\emptyset, \{b\}, \{a, b\},\{b, c\}, X\} \mbox{ and }$$
$$\mathcal{F}_{X/G}=\{\emptyset, \{[b]\}, \{[b], [c]\}, X/G\}.$$ Obviously, $F_{X/G}\in E(F_X)$.}
\end{example}

\begin{theorem}\label{theo1}
Let $L_0, L\in \mathcal{L}_{CD}$ and they have $DP$ and satisfy \emph{($\mathfrak{M}$)}. Then $L_0\in E(L)$
if and only if there exists a monotonic operator $G$ on $M(L)$ such that $L_0\cong F_{M(L)/G}$.
\end{theorem}
\proof
Notice that, from Lemma \ref{lemm 3}, we know that $L\cong F_{M(L)}$ and $L_0\cong F_{M(L_0)}$.

Now, suppose that $L_0\in E(L)$. Then there exists $\mathcal{S} \subseteq  \mathcal{F}_{M(L)}$ with $\emptyset, M(L) \in \mathcal{S}$ such that
$L_0\cong (\mathcal{S}, \supseteq)$ and $(\mathcal{S}, \supseteq)$ is a complete sublattice of $F_{M(L)}$.
Let $S=(\mathcal{S}, \supseteq)$. Then $L_0\cong S$. This means $(M(L_0), \leq)\cong (M(S), \supseteq)$.
Thus, $L_0\cong F_{M(S)}$ since $L_0\cong F_{M(L_0)}$. Therefore, it suffices to show that there exists a monotonic operator $G$ on $M(L)$ such that
\begin{equation}\label{he332}F_{M(S)}\cong F_{M(L)/G}.\end{equation}

Let $\mu: M(L)\rightarrow S$ be defined by \begin{equation}\label{he1}\mu(x)=\bigcap\{p\in \mathcal{S}\mid x\in p\}\end{equation} for all $x\in M(L)$.
Since $S$ is a complete sublattice of $F_{M(L)}$ and $M(L)\in \mathcal{S}$, we know that $\mu$ fulfills the conditions of Proposition \ref{pro 1}.
Thus $\mu$ is an $S$-fuzzy up-set on $M(L)$.

First, we shall prove that \begin{equation}\label{he2}M(S)=\mu(M(L)).\end{equation}

Let $x\in M(L)$. Note that $S$ has $DP$ since
$S\cong L_0$. Thus there exists a set $\{T_i\mid i\in I\}\subseteq M(S)$ such that $\mu(x)=\bigwedge_{i\in I}T_i$. Assume that $\mu(x)\notin M(S)$. Then for all $i\in I$, $T_i> \mu(x)$, i.e.,
\begin{equation}\label{he111}T_i\subsetneq \mu(x) \mbox{ for all }i\in I.\end{equation}
We claim that \begin{equation}\label{he110}I\neq \emptyset.\end{equation}
Indeed, if $I= \emptyset$ then $\mu(x)=\bigwedge \emptyset=1_S=\emptyset$ (in the complete sublattice $S=(\mathcal{S}, \supseteq)$, $\emptyset$ is the top element $1_S$), a contradiction since $x\in \mu(x)$ by formula (\ref{he1}).
Moreover, since $S$ is a complete sublattice of $F_{M(L)}$,
\begin{equation}\label{he112} \mu(x)=\bigwedge_{i\in I}T_i=\bigcup_{i\in I}T_i.\end{equation}
Thus, by formulas (\ref{he110}) and (\ref{he112}), it follows from $x\in \mu(x)$ that there exists $k\in I$ such that $x\in T_k$.
Further, by formula (\ref{he1}), we have $T_k\supseteq \mu(x)$, contrary to (\ref{he111}). Therefore, $\mu(x)\in M(S)$, i.e., $\mu(M(L))\subseteq M(S)$.

On the other hand, let $T\in M(S)$. Then $T\neq 1_S$ and $T\in \mathcal{S}$. Note that $0_S=M(L)$.
Thus, $T\subseteq M(L)$. Furthermore, by formula (\ref{he1}), $z\in \mu(z)\subseteq T$ for all $z\in T$.
Thus $T\subseteq \bigcup_{z\in T}\mu(z)\subseteq T$. Then $T=\bigcup_{z\in T}\mu(z)$, which means that $T=\bigwedge_{z\in T} \mu(z)$
since $S$ is a complete sublattice of $F_{M(L)}$.
As $T\in M(S)$ and $T\neq \emptyset$, there exists $y\in T$ such that $\mu(y)=T$. Hence, $T\in \mu(M(L))$. Therefore, $M(S)\subseteq \mu(M(L))$.

In view of the proof in the preceding two paragraphs, $M(S)=\mu(M(L))$.

Secondly, let \begin{equation}\label{che1}G: M(L)\rightarrow P_{M(L)} \mbox{ be defined by }G(x)=\mu(x)\end{equation} for all $x\in M(L)$. Then we shall prove \begin{equation}\label{he113}(M(S), \supseteq)\cong (G(M(L)), \supseteq).\end{equation}
Note that $S$ is a sub-poset of $P_{M(L)}$. As $\mu$
is an $S$-fuzzy up-set on $M(L)$, if $x\leq y$ then $G(x)=\mu(x) \leq \mu(y)=G(y)$, i.e.,
$G(x)\supseteq G(y)$. Thus, $G$ is a monotonic operator on $M(L)$ by Definition \ref{def 1}. From formula (\ref{che1}), $G(M(L))=\mu(M(L))$. Thus, by formula (\ref{he2}),
$M(S)=\mu(M(L))=G(M(L))$. Therefore, $(M(S), \supseteq)\cong (G(M(L)), \supseteq)$.

Finally, from Lemma \ref{lemma 1}, $(M(L)/G, \leq)\cong (G(M(L)), \supseteq)$. Thus, by formula (\ref{he113}), $(M(S), \supseteq)\cong (M(L)/G, \leq)$.
Therefore, $F_{M(L)/ G}\cong F_{M(S)}$, i.e., (\ref{he332}) is true.

Conversely, suppose that there exists a monotonic operator $G$ on $M(L)$ such that $L_0\cong F_{M(L)/G}$. Then by Theorem \ref{the0 6}, $F_{M(L)/G}\in E(F_{M(L)})$. Therefore, $L_0\in E(L)$ since $L\cong F_{M(L)}$.\endproof

Note that Example \ref{exa 1} also tell us that if $L_0$ is a finite sublattice of a distributive lattice $L$
then $L_0$ may be not in $E(L)$ generally. However, the following theorem will show us that all finite sublattices of
a complete atomic boolean lattice $L$ are in $E(L)$.

Let $A$ and $B$ be two sets. Then we denote $A\setminus B=\{x\in A\mid x\notin B\}$,
for convenience, if $B=\{b\}$ then we write $A\setminus B$ as $A\setminus b$.

\begin{theorem}\label{theo3}
Let $L$ be a complete atomic boolean lattice and let $L_0$ be a complete sublattice of $L$ and $L_0$ have $DP$ and satisfy \emph{($\mathfrak{M}$)}.
Then $L_0\in E(L)$.
\end{theorem}
\proof
Because $L$ is a complete atomic boolean lattice, we know that $L\cong P_{M(L)}$ and for all $a, b\in M(L)$,
\begin{equation}\label{he4}a \mbox{ and }b \mbox{ are incomparable. } \end{equation}
Thus $\mathcal{P}(M(L))=\mathcal{F}_{M(L)}$, and then \begin{equation}\label{he3}L\cong F_{M(L)}.\end{equation}
Furthermore, by Theorem \ref{theo 5}, $L$ has $DP$ and satisfies ($\mathfrak{M}$).
Therefore, by Theorem \ref{theo1}, we only need to construct a monotonic operator $G$ on $M(L)$ such
that \begin{equation}\label{eeel}L_0\cong F_{M(L)/G}.\end{equation}

First, as $L_0$ is a complete sublattice of $L$, it follows from formula (\ref{he3}) that there exists a sublattice $(\mathcal{S}, \supseteq)$ of ${F}_{M(L)}$ such that $L_0\cong (\mathcal{S}, \supseteq)$. Let $(\mathcal{S}, \supseteq)=S$, and then let $\mathcal{S}^{'}=\{T\setminus 1_S: T\in \mathcal{S}\}$
and $S^{'}=(\mathcal{S}^{'}, \supseteq)$. Clearly, $1_{S^{'}}=\emptyset$.
Moreover, let $f: S\rightarrow S^{'}$ be defined
by $$f(T)=T\setminus 1_S \mbox{ for all } T\in \mathcal{S}.$$
Obviously, $f$ is an isomorphism from $S$ to $S^{'}$, i.e., $S\cong S^{'}$.
Moreover, from the construction of $S^{'}$, we know that $\mathcal{S}^{'}\subseteq \mathcal{F}_{M(L)}$ since $\mathcal{F}_{M(L)}=\mathcal{P}(M(L))$.
Thus, as $S$ is a complete sublattice of $F_{M(L)}$, $S^{'}$ is also a complete sublattice of $F_{M(L)}$.

Secondly, since $0_{S^{'}}\neq \emptyset$ and $0_{S^{'}}\subseteq M(L)$, we know that $(0_{S^{'}}, \leq)$ is a sub-poset of $M(L)$.
Thus for all $x, y\in 0_{S^{'}}$, $x$ and $y$ are incomparable by (\ref{he4}). Let $\mu: (0_{S^{'}}, \leq)\rightarrow S^{'}$ be defined by
\begin{equation*}\mu(x)=\bigcap\{p\in \mathcal{S}^{'}\mid x\in p\}\end{equation*} for all $x\in 0_{S^{'}}$.
Note that $S^{'}$ is a complete sublattice of $F_{M(L)}$. Thus $\mu$ fulfills the conditions of Proposition \ref{pro 1}.
So that $\mu$ is an $S^{'}$-fuzzy up-set on $(0_{S^{'}}, \leq)$. Similar to the proof of (\ref{he2}), we can prove \begin{equation}\label{eeef}M(S^{'})=\mu(0_{S^{'}}).\end{equation}

Suppose that $w\in 0_{S^{'}}$ and let $G: M(L)\rightarrow P_{M(L)}$ be defined by
\begin{equation}\label{he5}
G(x)=\left\{
\begin{array}{rcl}
\mu(x)     &      & {x\in 0_{S^{'}},}\\
\mu(w)     &      & {x\in M(L)\setminus 0_{S^{'}}.}
\end{array}\right.
\end{equation}
Clearly, $G(M(L))=\mu(0_{S^{'}})$, which together with (\ref{eeef}) deduces that $G(M(L))=M(S^{'})$.
Furthermore, by formula (\ref{he4}), $G$ is a monotonic operator on $M(L)$.

Finally, from Lemma \ref{lemma 1}, $(M(L)/G, \leq)\cong (G(M(L)), \supseteq)$. Thus, $G(M(L))=M(S^{'})$
implies that $(M(S^{'}), \supseteq)\cong (M(L)/G, \leq)$. Hence, \begin{equation}\label{eeev}F_{M(L)/ G}\cong F_{M(S^{'})}.\end{equation}

On the other hand, from $L_0 \cong S$ and $S\cong S^{'}$, we know $L_0 \cong S^{'}$. Thus $(M(L_0), \leq)\cong (M(S^{'}), \supseteq)$. So that $F_{M(S^{'})}\cong F_{M(L_0)}$.
As $L_0$ has $DP$ and satisfies ($\mathfrak{M}$), we have that $L_0\cong F_{M(L_0)}$
by Lemma \ref{lemm 3}.
Therefore, by (\ref{eeev}), $F_{M(L)/ G}\cong F_{M(L_0)}\cong L_0$, i.e., (\ref{eeel}) holds.
\endproof

By Corollary \ref{cor4} and Theorem \ref{theo3}, the next corollary is obviously.

\begin{corollary}\label{cor6}
Let $L$ be a complete atomic boolean lattice. If $L_0$ is a finite sublattice of $L$ then $L_0\in E(L)$.
\end{corollary}

One can check that every sublattice $L_0$ of a finite
chain $L$ satisfies $L_0\in E(L)$. However, $L$ is not a complete atomic boolean lattice when $|L|\geq 3$.
\section{Conditions under which the poset of classes of monotonic operators is a lattice}

Let $\mathcal{H}(L)$ be the set of all monotonic operators on $M(L)$. For each $G\in \mathcal{H}(L)$,
we denote that $S_G=(\mathcal{S}_{G}, \supseteq)$ where $\mathcal{S}_{G}$ is defined by (\ref{equa 5}).
Let $L\in \mathcal{L}_{CD}$ and it have $DP$ and satisfy ($\mathfrak{M}$).

From Theorem \ref{theo1}, we know that the monotonic operator $G$ plays an important role in studying the structure of a lattice $L$ which can be represented as the collection of all up-sets. However, there may be two different monotonic operators $G_1, G_2\in \mathcal{H}(L)$
such that $X/G_1=X/G_2$. Then $S_{G_1}=S_{G_2}$ as illustrated by the following example.

\begin{example}\label{exa5}
\emph{Let us consider Example \ref{exam 2} again. Let $G_1: X\rightarrow P_X$ be a function
with $G_1(x)=\left\{
\begin{array}{rcl}
a_2     &      & {x=a,}\\
b_2     &      & {x=b,}\\
c_2     &      & {x=c.}
\end{array}\right.$ Obviously, $G_1$ is also a monotonic operator on $X$.
The poset $X/G_1$ is represented as Fig. 3.}
\par\noindent\vskip50pt
 \begin{minipage}{11pc}
\setlength{\unitlength}{0.75pt}\begin{picture}(600,120)
\put(270,110){\circle{4}}\put(273,115){\makebox(0,0)[l]{\footnotesize $[b]$}}
\put(270,75){\circle{4}}\put(252,80){\makebox(0,0)[l]{\footnotesize $[c]$}}
\put(270,40){\circle{4}}\put(273,45){\makebox(0,0)[l]{\footnotesize $[a]$}}
\put(270,42){\line(0,1){31}}
\put(270,77){\line(0,1){31}}
\put(256,20){$X/G_1$}
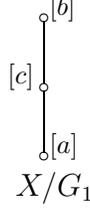
\put(171,0){\emph{Fig.3 Hasse diagram $X/G_1$}}
\end{picture}
\end{minipage}\\
\\
\emph{One can check that $X/G_1=X/G$. However, $G_1\neq G$.}
\end{example}

Therefore, the set of $\mathcal{H}(L)$ does not reflect the structure of $L$ really. Motivated by the forgoing reasons, in this section, we shall obtain an equivalence relation on $\mathcal{H}(L)$, show that classes of $\mathcal{H}(L)$ under an equivalence relation can be naturally ordered and give conditions under which the poset of classes of $\mathcal{H}(L)$ is a lattice.

\begin{definition}\label{defin 1}
\emph{Let $\approx$ be the relation on $\mathcal{H}(L)$, defined as:
$$G_1\approx G_2 \mbox{ if and only if }S_{G_1}=S_{G_2}.$$}
\end{definition}

From Definition \ref{defin 1}, we easily obtain the following lemma.
\begin{lemma}\label{lemma h1}
Let $L\in \mathcal{L}_{CD}$ and $L$ have $DP$ and satisfy \emph{($\mathfrak{M}$)}. Then, for any $G_1, G_2\in \mathcal{H}(L)$, \\
\emph{($i_1$)} $\approx$ is an equivalence relation on $\mathcal{H}(L)$.\\
\emph{($i_2$)} The set $\mathcal{H}(L)$ can be ordered: $[G_1]\leq [G_2]$ if and only if $\mathcal{S}_{G_1}\subseteq \mathcal{S}_{G_2}$
where $[G]=\{K\in \mathcal{H}(L)\mid K\approx G\}$.
\end{lemma}

\begin{theorem}\label{theorem 2}
Let $L$ be a finite distributive lattice. Then $(\mathcal{H}(L)/\approx, \leq)$ is a lattice.
\end{theorem}
\proof
By Corollary \ref{cor4} and Lemma \ref{lemm 3}, $L\cong F_{M(L)}$. Let $S=(\mathcal{S}, \supseteq)$ be a
sublattice of $F_{M(L)}$ with $\{\emptyset, M(L)\}\subseteq \mathcal{S}$. Then $S$ is a finite distributive lattice.
Again, by Corollary \ref{cor4} and Lemma \ref{lemm 3}, $S\cong F_{M(S)}$. Now, we shall prove that $(\mathcal{H}(L)/\approx, \leq)$ is a lattice
by two steps as below.

(I) There exists a monotonic operator $G\in \mathcal{H}(L)$ such that $S_G=S$, i.e., $(\mathcal{S}_G, \supseteq)=(\mathcal{S}, \supseteq)$.\\

Let $G$ be a function defined by (\ref{che1}) in which $\mu$ is given by (\ref{he1}). Then, by the proof of Theorem \ref{theo1},  we know that $G\in \mathcal{H}(L)$ and $\mu$ fulfills the conditions
of Proposition \ref{pro 1}. Let $p\in \mathcal{S}$. Then
$$p=\mu_p=\{x\in M(L)\mid \mu(x)\subseteq p\}.$$
Further, by (\ref{che1}), $p=\{x\in M(L)\mid G(x)\subseteq p\}$. Therefore,
\begin{equation}\label{che2}p=\{y\in[x]\mid x\in M(L), G(x)\subseteq p\}\end{equation} by the definition of $[x]$ in Lemma \ref{lemma 1}.

Let $T=\{[x]\mid x\in M(L), G(x)\subseteq p\}$. Suppose that $[y]\in M(L)/G$, $[x]\in T$ and $[y]\geq [x]$. Then $G(y)\subseteq G(x)$, obviously.
Note that $G(x)\subseteq p$. Thus $G(y)\subseteq p$, which means that $[y]\in T$.
Then $T$ is an up-set of $M(L)/G$. Hence, from formulas (\ref{equa 4}) and (\ref{che2}), we have $p=S_T$.
Thus $p\in \mathcal{S}_G$ by (\ref{equa 5}). Therefore, \begin{equation}\label{FFE}\mathcal{S}\subseteq \mathcal{S}_G.\end{equation}

On the other hand, by the proof of Theorem \ref{theo1}, $F_{M(S)}\cong F_{M(L)/G}$.
Moreover, by (\ref{equa 7}), $F_{M(L)/G}\cong S_G$. Thus $F_{M(S)}\cong S_G$, which together
with $S\cong F_{M(S)}$ yields that $S\cong S_G$. Therefore, by (\ref{FFE}), $\mathcal{S}=\mathcal{S}_G$, i.e., $S=S_G$.

(II) $(\mathcal{H}(L)/\approx, \leq)$ is a lattice.\\

The proof of (II) is made in two steps.

(i) By (I), there exists $U, V\in \mathcal{H}(L)$ such that $S_U=F_{M(L)}$ and $S_V=(\{\emptyset, M(L)\}, \supseteq)$, respectively. Obviously, $[U]$ and $[V]$ are the top and bottom elements of
$(\mathcal{H}(L)/\approx, \leq)$, respectively.

(ii) Suppose that $[G_1], [G_2] \in \mathcal{H}(L)/\approx$. From the proof of Theorem \ref{the0 6}, both $S_{G_1}$ and $S_{G_2}$ are the sublattices of $F_{M(L)}$ and
$\{\emptyset, M(L)\}\subseteq \mathcal{S}_{G_1}\bigcap\mathcal{S}_{G_2}$. So that $(\mathcal{S}_{G_1}\bigcap\mathcal{S}_{G_2}, \supseteq)$ is a sublattice of $F_{M(L)}$. Further, by (I), there exists $G\in \mathcal{H}(L)$ such that
$S_G=(\mathcal{S}_{G_1}\bigcap\mathcal{S}_{G_2}, \supseteq)$. Therefore, $[G]=[G_1]\wedge [G_2]\in \mathcal{H}(L)/\approx$ by Lemma \ref{lemma h1}.

By (i) and (ii), we conclude that $(\mathcal{H}(L)/\approx, \leq)$ is a lattice.
\endproof

\section{Conclusions}
As it is well known, the famous problem of Dilworth says that for an algebraic distributive lattice $L$,
whether there exists a lattice $L_0$ such that $L\cong Con(L_0)$, often referred to as CLP.
In 2007, Wehrung constructed an infinite algebraic distributive lattice $L$ satisfying that $L\ncong Con(L_0)$ for all
lattice $L_0$ (see \cite{Wehrung}).

On the other hand, Gr\"{a}tzer proved that every finite distributive lattice $L$ can be represented as the congruence
lattice of a finite semimodular lattice (even a finite rectangular lattice) $L_0$ with $(M(L), \leq)\cong (M(Con(L_0)), \leq)$ (see \cite{Gratze}).
In this paper, we have proved that if $L_1, L_2\in \mathcal{L}_{CD}$ and they have $DP$ and satisfy ($\mathfrak{M}$), then $L_1 \cong L_2$
if and only if $(M(L_1), \leq)\cong (M(L_2), \leq)$ (Corollary \ref{cor3}). Therefore, a naturally problem is: Is there a
semimodular lattice (even a rectangular lattice) $L_0$ such that $L\cong Con(L_0)$ if $L\in \mathcal{L}_{CD}$ and $L$
has $DP$ and satisfies ($\mathfrak{M}$)?

\end{document}